\documentclass{article}
\usepackage{verbatim}
\usepackage{booktabs}
\usepackage{hyperref}
\usepackage{mathrsfs}
\usepackage{tikz}
\usepackage{amsmath,amsthm,amssymb}
\usepackage{cleveref}
\usepackage{enumerate}

\newtheorem{theorem}{Theorem}
\newtheorem*{theorem*}{Theorem}

\newtheorem{lemma}[theorem]{Lemma}

\theoremstyle{definition}
\newtheorem{definition}[theorem]{Definition}

\DeclareMathOperator{\ob}{ob}

\DeclareMathOperator{\avoidOnly}{avoidOnly}
\DeclareMathOperator{\avoidNone}{avoidNone}
\DeclareMathOperator{\noObligations}{noObligations}
\newcommand{\fivecweak}{
	If $Y\in\ob(X)$
	and $Z\in\ob(X)$
	and $X\cap Y\cap Z\ne\emptyset$,
	then $Y\cap Z\in\ob(X)$.
}

\newcommand{\fivef}{
	If $X\in\ob(Y)$ and $X\in\ob(Z)$ then $X\in\ob(Y\cup Z)$.
}
\newcommand{\fiveg}{
	If $Y\in\ob(X)$ and $Z\in\ob(Y)$ and $X\cap Y\cap Z\ne\emptyset$, then $Y\cap Z\in\ob(X)$. (A form of transitivity for obligations.)
}

\newcommand{\fivee}{
	If $Y\subseteq X$ and $Z\in\ob(X)$ and $Y\cap Z\ne\emptyset$, then $Z\in\ob(Y)$.
}
\newcommand{\fived}{
	If $Y\subseteq X$ and $Y\in\ob(X)$ and $X\subseteq Z$, then $(Z\setminus X)\cup Y\in\ob(Z)$.
}
\newcommand{\fivec}{
	If $Y\in\ob(X)$ and $Z\in\ob(X)$ then $Y\cap Z\in\ob(X)$.
}
\newcommand{\fiveb}{
	If $Y\cap X=Z\cap X$ then $Y\in\ob(X)$ iff $Z\in\ob(X)$.
}
\newcommand{\fivea}{
	$\emptyset\not\in\ob(X)$.
}

\title{Classification and deontic explosion for contrary-to-duty obligations}
\author{Bjørn Kjos-Hanssen \\ University of Hawai\textquoteleft i at M\=anoa \\ bjoernkh@hawaii.edu}
\author{Bjørn Kjos-Hanssen}
\date{}

\begin{document}
\maketitle

\begin{abstract}
Carmo and Jones have presented a sequence of candidate axiom systems for conditional obligation between 1997 and 2022.
For their most recent system we demonstrate a limited form of deontic explosion: given that a student does not get the highest possible grade on a test, any other passing grade is acceptable.

In addition to that negative result, we give a positive one: revisiting the strongest version of Carmo and Jones' 1997 system, we provide a surprising classification of all satisfying models in terms of a single forbidden possible world.
\end{abstract}
\section{Introduction}
Carmo and Jones in 2022 \cite{MR4500520} proposed certain axioms 5(a)--(g) for a relation $\ob$ that holds between sets of possible worlds $X$ and $Y$ if $X$ is obligatory in the context $Y$. It was the latest iteration in a sequence of systems \cite{CJ96,CJ02,MR3063042}.

We will exhibit an anomaly therein: a weak form of conditional deontic explosion. Given that something is somewhat desirable (passing a course with a grade of C, say) and given that the most desirable outcome (the grade of A) is unavailable, the somewhat desirable outcome becomes obligatory.

We then show that despite this explosion, the systems of Carmo and Jones have interesting mathematical content. For the strongest system we provide a full classification of its models; for weaker versions we characterize the least models (under inclusion) satisfying the axioms and basic contrary-to-duty assumptions.

\section{The anomaly}
Let $W$ be a finite set of possible worlds of a given model, and let $\mathscr P$ denote the power set operation.

Combining the systems from several papers, the full list of Carmo and Jones'
conditions
on a function $ob : \mathscr P(W) \to \mathscr P(\mathscr P(W))$ is as follows. 
\begin{itemize}
	\item Axiom 5(a): \fivea (Obligations cannot be impossible to fulfil.)
	\item Axiom 5(b): \fiveb
		(Whether $Z$ is obligatory in context $X$ depends only on $Z$ through $Z\cap X$.)
	\item Axiom 5(c$^-$) (\cite[page 319]{CJ02};
		finite version of 5(c$^*$) in \cite{MR3063042}, which is called 5(c) in \cite{MR4500520}):
	\fivecweak (Closure under intersection.)
	\item Axiom 5(c) (\cite{CJ96} and \cite[page 287]{CJ02}): \fivec (Closure under intersection, even for disjoint sets.)
	\item Axiom 5(d): \fived
		(If $Y$ is obligatory in context $X$ then a form of the material conditional ``$X \to Y$'' holds in a wider context.)
	\item Axiom 5(e): \fivee
		(If $Z$ is obligatory in the context $X$ then it remains obligatory in any subcontext in which it is still possible.)
	\item Axiom 5(f): \fivef
		(If $X$ is obligatory in both of the contexts $Y$ and $Z$, then it is obligatory in their union.)
	\item Axiom 5(g): \fiveg
\end{itemize}

\paragraph{The argument.} Axioms 5(b)(e)(f) will now be used to derive an anomaly.

Suppose that James is taking an exam on which the possible grades, in decreasing order of quality, are A, B, C, D and F.

If seems reasonable to assume that given that James' grade is A or B, it ought to be A.

Moreover, given that the grade is C or D, it ought to be C.

Then by axiom 5(b), given that the grade is C or D, it ought to be A or C.
Moreover, also by 5(b), given that the grade is A or B, it ought to be A or C.
In other words, James ought to make sure the proposition ``the grade is A or C'' is true (of course, ideally by getting an A).

It also follows that the grade ought to be A or C given that it is A, B, C, or D,
this time using axiom 5(f) applied to the two previous statements.
(If this is getting hard to follow, the reader may rest assured that a formalization is available, see \Cref{sa}.)

So far this is not entirely unintuitive. But now we use 5(e), and conclude that
since A or C was obligatory in the context A, B, C, or D
it must remain obligatory in a more restrictive context in which it is still possible, namely B, C, or D.

Finally, using 5(b) again, we conclude that given that the grade is B, C, or D, it ought to be C.
While there is such a concept as ``gentleman's C'',
surely this is a paradox.\footnote{Traditionally a \emph{Gentleman's C}
is given to someone who deserves D or F but for extraneous reasons is deemed worthy of a C.
Here this is reversed, as a C is deemed obligatory even when a B is available.
}

\paragraph{A literary perspective.} 
Let us consider the anomaly above in another context, beside students and their grades.

Suppose James has the following options.
\begin{enumerate}[A]
\item Marriage to Alice, his favorite.
\item Marriage to Alice's sister Beatrice, or one of several other women he also quite likes.
\item Remaining a lonely bachelor.
\item Marriage to Deirdre, whom he despises.
\end{enumerate}

Then our conclusion is like that of the Bee Gees in 1977:
\begin{quote}
If I can't have you, I don't want nobody.
\end{quote}

\paragraph{From despondency to mathematics.}
Hage \cite{Hage2000-HAGDNE} thinks the task Carmo and Jones and others have set themselves is impossible.

From Hage's book review, in the context of discussing a man who
(i) should visit his neighbor, (ii) if he does visit, should call to say that he is coming, and
(iii) is not in fact coming:

\begin{quote}
	[...] Actually, A does not assist his neighbours and
	therefore should not call them. The two intuitions presuppose a different role for
	deontic logic, namely, reasoning about what is ideally the case, and reasoning about
	what ought to be done in the real world. These two roles are hard to reconcile in
	one logic, unless the logic combines the two kinds of reasoning in different parts.
	The paper by Carmo and Jones in the present volume provides such a combination
	logic.
	Nevertheless, many have attempted to do what is, in my opinion, \textbf{undoable}\footnote{[emphasis ours]} and
	this has lead to many modifications of the Standard System of Deontic Logic
	(Hilpinen, 1971, p. 13f.; Chellas, 1980, p. 190f.). These new systems have lead
	to new paradoxes that deal with CTD obligations.
\end{quote}

Given Hage's sentiment
and the anomaly we have presented, one might temporarily become despondent.
If deontic systems always have flaws, why pursue them?
However, in the course of uncovering this paradox I also discovered interesting \emph{mathematical structure} in the axioms 5(a)--5(g).

\begin{definition}
Let us say that a function $\ob$ is the \emph{least} member of a family $\mathcal O$ if $\ob\in\mathcal O$,
and for each $\ob'\in\mathcal O$ and each $X$ and $Y$,
if $X\in\ob(Y)$ then $X\in\ob'(Y)$.
\end{definition}

The least model of 5(b) given certain ``oughts'' has a very natural characterization,
as does the least model of 5(b), 5(d), and 5(f).

To be specific, let us write $\mathrm{Ought}(A\mid B)$ to mean that for each $X\subseteq B$,
if $A\cap X\ne \emptyset$ then $A\in\ob(X)$.
This is the semantic condition used by Carmo and Jones for the conditional obligation operator $O(A\mid B)$.
A pair of oughts $(\mathrm{Ought}(A\mid W)$, $\mathrm{Ought}(B\mid W \setminus A))$
forms a basic contrary-to-duty obligation of $B$ given that our duty $A$ has failed to be observed.

\begin{definition}\label{canon2}
The set canon$_2$ A B X consists of all contexts that are obligatory in the model canon$_2$ A B.

It is defined to be: if $X \cap B = \emptyset$ then $\emptyset$, else:
if $X \cap A = \emptyset$ then $\{T \mid X \cap B \subseteq T\}$, else $\{T \mid X \cap A \subseteq T\}$.
\end{definition}

\begin{definition}\label{canon2_II}
The set canon$_{2,II}$ A B X consists of all contexts that are obligatory in the model canon$_{2,II}$ A B.

It is defined to be: if $X \cap B = \emptyset$ then $\emptyset$, else:
if $X \cap A = \emptyset$ then $\{Y \mid X \cap B = X \cap Y\}$, else $\{Y \mid X \cap A = X \cap Y\}$.
\end{definition}
In other words, if $Y \in\ob(X)$ under a canon$_{2,II}$ model
then generically $X\cap Y$ consists exactly of the most desirable worlds in $X$.

\begin{theorem}\label{characterize_canon2_II}
The least model $\ob$ of $\mathrm{Ought}(A\mid W)$, $\mathrm{Ought}(B\mid W \setminus A)$, and axiom 5(b) is canon$_{2,II}$ A B.
\end{theorem}

Even though the model arises from assuming 5(b) only, it also satisfies axioms 5(a), 5(c), 5(e) and 5(g).
This indicates a certain robustness of our definitions.

\begin{theorem}\label{characterize_canon2}
The least model of axioms 5(b), 5(d), 5(f) and the two ``oughts'' in
\Cref{characterize_canon2_II}
is canon$_2$ A B.
\end{theorem}

Since axiom 5(f) follows from 5(a)(b)(c)(d), as shown by Carmo and Jones,
the model can alternatively be characterized as the least one satisfying the latter four axioms and the two specified Oughts.

The two families of models in \Cref{characterize_canon2_II} and \Cref{characterize_canon2}
represent two alternative approaches to contrary-to-duty obligations as discussed in \cite{MR3607634}
(called II and I there, respectively).
They do not exhaust the interesting models by any means: for instance, we may have conflicting obligations.
This may lead us to prefer the weak version of 5(c) to the strong version.
For a concrete example, suppose that James has received acceptances on separate marriage proposals to both Alice and Beatrice,
but cannot marry them both.

\paragraph{Classification.}
In mathematics, classification theorems are fairly common.
For example, the finite abelian groups have a straightforward characterization and the finite simple groups a complicated one.
Vector spaces over a fixed field $\mathbb F$ of finite dimension are characterized by their dimension $d$,
hence the single parameter $d$ determines the space up to isomorphism.

In deontic logic, where axioms and rules are added based on moral intuitions,
we should perhaps not expect structures that are mathematically natural enough to be classifiable.

However, for the full theory of Carmo and Jones's axioms 5(a)--5(e),
with the strong version of 5(c) in which we do not impose non-disjointness,
we can characterize its models completely.
The only nontrivial ones basically say that there is just one bad world and the only obligation is to avoid it:
\begin{theorem}\label{models_ofCJ_1997_equiv}
	Let $W$ be a finite set of possible worlds and let $\ob : \mathscr P (W)\to \mathscr P(\mathscr P(W))$.
	The following are equivalent:
	\begin{enumerate}
		\item $\ob$ satisfies the full system suggested in \cite{CJ96}: 5(a), (b), (c), (d) and (e).
		\item One of the following three holds:
		\begin{enumerate}
			\item $\ob(X)=\emptyset$ for all $X$.
			\item $\ob(X)=\{Y \mid \emptyset \ne X \subseteq Y\}$ for all $X$.
			\item There is a distinguished possible world $a$ (the ``forbidden'' world) such that for all $X$,
			$\ob(X)=\{Y \mid X \cap Y \ne \emptyset \text{ and } X\setminus \{a\}\subseteq Y\}$.
		\end{enumerate}
	\end{enumerate}
\end{theorem}
Details on the proof of Theorem \ref{models_ofCJ_1997_equiv} can be found in the Appendix.
Finally, we note that the semantic conditions above are independent, as shown by the following result (proved in \cite{deontic}).
\begin{theorem}
	The axioms 5(a), 5(b), 5(c) (strong version), 5(d), and 5(e) are independent
	in the sense that no single axiom follows from the other four axioms.
\end{theorem}

\section{Conclusion}
We have seen that the system CJ97 has a mathematical model classification,
and that the weaker subsystem using only axioms 5(b)(e)(f) already suffers from a form of deontic explosion.
Future work that may be undertaken includes:
\begin{itemize}
	\item First and foremost, to continue the exploration of subsystems of CJ97 or similar systems,
		and assessing their adequacy for deontic logic.
	\item Second, given that the models of the system CJ97 have been classified exactly,
		it would be interesting to see if a subsystem,
		such as the one obtained by replacing the strong version of axiom 5(c) with the weak version, also has a classification.
\end{itemize}

\section{Acknowledgments}\label{sa}

All mathematical claims above are verified in the proof assistant \textsc{Lean}, see \cite{deontic}. The main argument for our anomaly was discovered by carefully analyzing some output from a script in the computer mathematics system \textsc{Maple} by Maplesoft \cite[Appendix]{K96}.

This work was partially supported by a grant from the Simons Foundation (\#704836 to Bj{\o}rn Kjos-Hanssen).

\bibliographystyle{alpha}
\bibliography{deontic}

\appendix
\section{Details of the characterization of CJ97}
By CJ97 we mean the system consisting of axioms 5(a)(b)(c)(d)(e), suggested in \cite{CJ96},
in particular using the strong version of 5(c).

\begin{definition}\label{CX}
	Conditional explosion for ob is the statement that for all subsets $A,B,C$ of the set of possible worlds $W$,
	if $A \in \ob(C)$ and $B \cap A^c \cap C \ne \emptyset$ then $B \in \ob(A^c \cap C)$.
\end{definition}

\begin{theorem}\label{conditional_explosion}
	If ob satisfies axioms 5(a)(b)(d)(e) then ob satisfies conditional explosion.
\end{theorem}

\begin{lemma}\label{single_ob_pair}
	Assume 5(d), 5(e), and 5(c$^-$).
	If $a_1\ne a_2$, $a_1\notin A$, $a_2\notin A$, $\{a_1,a_2\}\in\ob(\{a_1,a_2\})$,
	and $A\in\ob(X)$ whenever $A\subseteq X$, then
		$\{a_1\} \in \ob(\{a_1, a_2\})$ and $\{a_2\} \in \ob(\{a_1, a_2\})$.
\end{lemma}

\begin{lemma}[proof uses \Cref{single_ob_pair}]\label{semiglobal_holds}
	Assume 5(a)(c)(d)(e).
	If $a_1\ne a_2$, $a_1\notin A$, $a_2\notin A$, $\{a_1,a_2\}\in\ob(\{a_1,a_2\})$,
	then it cannot be that $A\in\ob(X)$ whenever $A\subseteq X$.
\end{lemma}

\begin{lemma}[proof uses \Cref{conditional_explosion} and \Cref{semiglobal_holds}]\label{global_holds_specific}
	Assume all axioms 5(a)(b)(c)(d)(e).
	If  $a_1 \notin A$ and $a_2 \notin A$ and $a_1 \ne a_2$ then it cannot be that $A\in\ob(X)$ whenever $A\subseteq X$.
\end{lemma}

\begin{definition}\label{bad}
A world $a$ is bad if  $\exists X, a \in X \wedge X \setminus \{a\} \in \ob(X)$.
\end{definition}
\begin{definition}\label{quasibad}
The world $a$ is quasibad if $\exists X, \exists Y, a \in X \setminus Y \wedge Y \in \ob(X)$.
\end{definition}

Thus, a world $a$ is bad if in some context there is an obligation to simply avoid $a$.
For example, if there is an obligation ``do not go to war'' then the world representing ``going to war with Syria''
is \emph{quasibad}, but it is not \emph{bad} unless there is also the specific obligation ``do not go to war with Syria''.
In ``reasonable'' systems this distinction would perhaps not need to be made,
but here we are in the process of proving that a certain system CJ97,
comprising all of 5(a)(b)(c)(d)(e), is not reasonable.

\begin{lemma}\label{obSelfSdiff_of_bad}
	Assume 5(b)(d)(e).
	If $a$ is bad and $Y \setminus \{a\} \ne \emptyset$ then
	$Y \setminus \{a\} \in \ob(Y)$.
	\end{lemma}

Lemma \ref{obSelfSdiff_of_bad} says that badness of the world does not depend on context.

\begin{lemma}\label{obSelf_of_obSelf}
	Assume 5(b)(d)(e).
	If  $X \in \ob(X)$ and $Y \ne \emptyset$ then $Y \in \ob(Y)$.
\end{lemma}

Lemma \ref{obSelf_of_obSelf} says that if any context is obligatory relative to itself, then they all are.

\begin{lemma}[proof uses \Cref{obSelf_of_obSelf}]\label{obSelf_of_obSelfSdiff}
	Assume 5(b)(d)(e).
	If $\emptyset \ne X \setminus \{a\} \in \ob(X)$ then $X \in \ob(X)$.
\end{lemma}

Lemma \ref{obSelf_of_obSelfSdiff} says that if there is a bad world, then the corresponding context is self-obligatory.

\begin{lemma}\label{obSelf_univ}
	Assume 5(a)(d)(e).
	If $W \setminus \{a\} \in \ob(W)$ then $W \in \ob(W)$.
\end{lemma}

\begin{lemma}[proof uses \Cref{obSelf_univ}]\label{obSelf_of_bad.single}
Assume 5(a)(b)(d)(e).
If $a$ is bad then $\{a\} \in \ob(\{a\})$.
\end{lemma}

Lemma \ref{obSelf_univ} is technicality: if $a$ is bad in the global context then the global context is self-obligatory.
Lemma \ref{obSelf_of_bad.single} is another technicality: even if $a$ is bad, it is still self-obligatory.

We say that a set $A$ is a subsingleton if its cardinality is at most 1,
and a cosubsingleton if $W \setminus A$ is a subsingleton.

\begin{lemma}\label{local_of_global}
	Assume 5(a)(d)(e).
	Suppose that for all contexts $A$, if $A$ is obligatory in all larger contexts $X\supseteq A$, then $A$ is a cosubsingleton.
	Then for all $B$ and $C$, if $B\subseteq C$ is obligatory relative to $C$ then $C\setminus B$ is a subsingleton.
\end{lemma}

Lemma \ref{local_of_global} is a ``global-to-local'' principle,
allowing us to conclude a fact about an arbitrary context $C$ from a fact about the global context.

The antecedent of Lemma \ref{local_of_global} is provided by Lemma \ref{global_holds}
and hence the consequent is provided by Lemma \ref{local_holds}.
\begin{lemma}[proof uses \Cref{global_holds_specific}]\label{global_holds}
	Assume all axioms 5(a)(b)(c)(d)(e).
	For all contexts $A$, if $A$ is obligatory in all larger contexts $X\supseteq A$, then $A$ is a cosubsingleton.
\end{lemma}

\begin{lemma}[proof uses \Cref{global_holds,local_of_global}]\label{local_holds}
	Assume all axioms 5(a)(b)(c)(d)(e).
	For all $B$ and $C$, if $B\subseteq C$ is obligatory relative to $C$ then $C\setminus B$ is a subsingleton.
\end{lemma}

We now come to define the three kinds of models, i.e., possible functions ob, that the system CJ97 has.

\begin{definition}\label{avoidOnly}
In the model \emph{avoidOnly} with distinguished world $e$, the contexts that are obligatory relative to $X$ are defined by
\[
\avoidOnly_e(X) = \{Y  | X \cap Y \ne \emptyset \wedge X \setminus \{e\} \subseteq  X \cap Y\}.
\]
\end{definition}

\begin{definition}\label{avoidNone}
In the model \emph{avoidNone}, the contexts that are obligatory relative to $X$ are defined by
\[
	\avoidNone(X) = \{Y \mid X \ne \emptyset \wedge Y \supseteq X\}.
	\]
\end{definition}

\begin{definition}\label{noObligations}
In the model \emph{noObligations} there are no obligations at all:
\[
	\noObligations(X) = \emptyset.
\]
\end{definition}

We can view \emph{avoidNone} as a special case of \emph{avoidOnly} where the world to be avoided does not exist in the model.

We prove several technical lemmas, culminating in \Cref{models_ofCJ_1997}:
\begin{lemma}[proof uses \Cref{obSelf_of_obSelf,local_holds}]\label{avoidNone_of_no_quasibad}
	Assume all axioms 5(a)(b)(c)(d)(e).
	If $\ob \ne \noObligations$ and there are no quasibad worlds, then $\ob = \avoidNone$.
\end{lemma}

\begin{lemma}[proof uses \Cref{local_holds}]\label{unique_bad}
	Assume all axioms 5(a)(b)(c)(d)(e).   
	If $a$ and $b$ are bad then $a = b$.
\end{lemma}

\begin{lemma}[proof uses \Cref{unique_bad}]\label{bad_cosubsingleton_of_ob}
	Assume 5(a)(b)(d)(e).
	If $a$ is bad and $X \cap Y \in \ob(X)$ then
		$X \cap Y = X$ or $X \cap Y = X \setminus \{a\}$.
\end{lemma}

\begin{lemma}[proof uses \Cref{obSelf_of_obSelfSdiff,obSelfSdiff_of_bad}]\label{obSelf_of_bad.nonsingle}
	Assume 5(b)(d)(e).
If $a$ is bad, $X \ne \{a\}$ and $X \ne \emptyset$ then $X \in \ob(X)$.
\end{lemma}

\begin{lemma}[proof uses \Cref{obSelf_of_bad.nonsingle,obSelf_of_bad.single}]\label{obSelf_of_bad}
	Assume 5(a)(b)(d)(e).
	If there exists a bad world at all, and $X \ne \emptyset$, then $X \in \ob(X)$.
\end{lemma}

\begin{lemma}[proof uses \Cref{bad_cosubsingleton_of_ob}]\label{sub_avoidOnly_of_bad}
	Assume all axioms 5(a)(b)(c)(d)(e).
	If $a$ is bad then for all $Y$, $\ob(Y) \subseteq \avoidOnly_a(Y)$.
\end{lemma}

\begin{lemma}[proof uses \Cref{obSelf_of_bad}]\label{avoidOnly_sub_of_bad}
	Assume 5(a)(b)(d)(e).
	If $a$ is bad then for all $Y$, $\avoidOnly_a(Y) \subseteq \ob(Y)$.
\end{lemma}
Unlike Lemma \ref{sub_avoidOnly_of_bad}, Lemma \ref{avoidOnly_sub_of_bad} does not require any form of axiom 5(c).

\begin{lemma}[proof uses \Cref{sub_avoidOnly_of_bad,avoidOnly_sub_of_bad}]\label{avoidOnly_of_bad}
	Assume all axioms 5(a)(b)(c)(d)(e).
	If $a$ is bad then $\ob = \avoidOnly_a$.
\end{lemma}

\begin{lemma}\label{bad_of_quasibad}
	Assume axioms 5(a)(b)(e)\footnote{The author's original proof of this lemma used all axioms 5(a)(b)(c)(d)(e).
	An anonymous referee proved that the reduced axiom set suffices.}.
	If a world $a$ is quasibad then $a$ is bad.
\end{lemma}

Finally, we obtain \Cref{models_ofCJ_1997_equiv}:
\begin{theorem*}[proof uses \Cref{avoidNone_of_no_quasibad,avoidOnly_of_bad,bad_of_quasibad}]\label{models_ofCJ_1997}
	Every model of axioms 5(a)(b)(c)(d)(e) is either $\avoidOnly_a$ for some bad world $a$,
	$\avoidNone$, or $\noObligations$.
\end{theorem*}

\end{document}